\documentclass{article}
\usepackage[utf8]{inputenc}
 \usepackage{arxiv}

\usepackage{amsmath}
\usepackage{amsfonts}
\usepackage{amssymb}
\usepackage{amsthm}
\usepackage{latexsym}
\usepackage{amscd}
\usepackage{hyperref}
\usepackage{tikz-cd}
\usetikzlibrary{decorations.pathmorphing}
\usetikzlibrary{decorations.pathreplacing}
\usepackage{url}
\usepackage{float}
\usepackage[T1]{fontenc}

\allowdisplaybreaks
\newtheorem{thm}{Theorem}[section]
\newtheorem{prop}{Proposition}[section]

\theoremstyle{definition}
\newtheorem{defi}{Definition}[section]

\title{Boundary value problem of magnetically insulated diode: existence of solutions and complex bifurcation}

\newif\ifuniqueAffiliation
\uniqueAffiliationtrue

\ifuniqueAffiliation % Standard variant of author block

\title{Boundary value problem of magnetically insulated diode: existence of solutions and complex bifurcation}

\renewcommand{\shorttitle}{BVP of MID: Solution and Bifurcation}
\newif\ifuniqueAffiliation
\uniqueAffiliationtrue

\ifuniqueAffiliation % Standard variant of author block
\author{Denis N.~Sidorov \\
	Department of Applied Mathematics\\
	MESI SB Russian Academy of Sciences\\
	Irkutsk, Russia \\
	\texttt{dsidorov@isem.irk.ru} \\
	\AND 
	Aleksander V.~Sinitsyn \\
	Departamento de Matem\'aticas\\
	Universidad Nacional de Colombia\\
	Bogot\'a, Columbia \\
	\texttt{asinitsyne@unal.edu.co} \\
	 \AND
	Toledo Leguizam\'on \\
Ingeniería de Sistemas y Computación\\
	Universidad Nacional de Colombia\\
	Bogot\'a, Columbia \\
	\texttt{otoledo@unal.edu.co} \\
	\AND
	Liguo Wang \\
Department of Electrical Engineering and Automation\\	
Harbin Institute of Technology \\
Harbin, PRC\\
\texttt{wlg2001@hit.edu.cn}
}

\begin{document}
\maketitle

%abstract
\begin{abstract}
The paper focuses on the stationary self-consistent problem of magnetic insulation for a vacuum diode with space-charge limitation, described by a singularly perturbed Vlasov-Maxwell system of dimension 1.5. The case of insulated diode when the electrons are deflected back towards the cathode at the point $x^{*}$ is considered. First,  the  initial VM system is reduced to the nonlinear singular limit system of ODEs for the potentials of electric and magnetic fields. The second step deals with the limit system's reduction to the new nonlinear singular ODE equation for effective potential $\theta(x)$. 
The existence of non-negative solutions is proved for the last equation on the interval $[0, x^{*})$ where $\theta(x)>0$. The most interesting and unexplored case is when $\theta(x)<0$ on the interval $(x^{*}, 1]$ and corresponds to the case of an insulated diode.
For the first time, a numerical analysis of complex bifurcation of solutions in insulated diode is considered for $\theta(x)<0$ depending on parameters and boundary conditions. Bifurcation diagrams of the dependence of solution $\theta(x)$ on a free point (free boundary) $x^{*}$ were constructed. Insulated diode spacing is found.\end{abstract}

%keywords
\keywords{relativistic Vlasov-Maxwell system, magnetic insulation, effective potential, insulated diode, initial value problem, singular boundary value problem, contractive mapping, fixed point theorem, complex numerical bifurcation.}
%Acknowledgments

%for citation

\section{Motivation}
Power electronics is a cornerstone of contemporary power systems, bringing theoretical physics and practical engineering to enable efficient energy conversion and transmission. Among its critical challenges is ensuring
the stability and reliability of high-voltage devices, such as thouse employed in high-voltage direct current (HVDC) systems. A fascinating intersection of applied mathematics and plasma physics arises in the study of magnetic insulation in vacuum diodes -- a phenomenon where magnetic fields surpass electron backflow, 
allowing these devices to operate efficiently at extreme voltages. The seminal work of Langmir and Compton \cite{1} provides the fundamental theory to model electron behavior.   

The Vlasov-Maxwell theory provides a rigorous kinetic framework to model the collective behavior of 
charged particles (electrons) and electromagnetic fields. Here readers may refer to  part 3 of  book \cite{sidsidsin}, \cite{8}
and references therein including works \cite{1,3,7}. The generic theory of the parametric families of small branching (bifurcation) solutions of nonlinear differential equations is proposed in \cite{1, sidsid} and applied to
magnetic insulation problem. The reduction of the boundary value problem for noninsulated magnetic regime in a vacuum diode to a singular system of nonlinear Fredholm equations is fulfilled in \cite{sym20}.

The nonlinear dynamics of these models often leads to bifurcations -- sudden qualitative shifts in system behavior caused by small parameter changes. Bifurcations analysis  is crucial for predicting instabilities, optimizing insulation thresholds, and preventing catasctrophic failures in devices like vacuum diodes. 
It is to be noted that electron transport in high energy devices such as vacuum diodes exhibits many nonlinear phenomena due to the extremely high applied voltages. 
One of these effects is the saturation of the current due to the self-consistent electric and magnetic field. 

The effect of magnetic insulation consists in that the electrons emitted from cathode cannot reach the anode due to the extremely high applied electric and magnetic field; they are reflected by the magnetic forces back to the cathode. Here two basic regimes are possible: the first, when electrons reach the anode ("noninsulated" diode) and the second one, when electrons rotate back to the cathode ("insulated" diode).

\section{Setting of the problem}
We consider a plane diode consisting of two perfectly conducting electrodes, a cathode $(X=0)$ and anode $(X=L)$. The system is described by the 1.5 dimensional VM model:

\begin{align*}
\,&V_{X}\frac{\partial F}{\partial X}+e\Biggl(\frac{d\Phi}{d X}-V_{Y}\frac{d A}{d X}\Biggr)\frac{\partial F}{\partial P_{X}}+eV_{X}\frac{d A}{d X}\frac{\partial F}{\partial P_{Y}}=0, \\
\,&\frac{d^{2}\Phi}{d X^{2}}=\frac{e}{\epsilon_{0}}N(X), \qquad X\in (0,L), \\
\,&\frac{d^{2}A}{d X^{2}}=-\mu_{0}J_{Y}(X), \qquad X\in(0,L).
\end{align*}

% insert 1

Here, $F(X,P_{X}, P_{Y})$ is electron distribution function, $\Phi,$ $A$ are potentials of electromagnetic field, $\varepsilon_0,$ $\mu_0$  are the vacuum  permittivity and permeability.

After appropriate scaling and taking the limit $\varepsilon\to0$, we obtain the limit system:

\begin{equation}
\begin{array}{l}
\frac{d^{2}\varphi}{dx^{2}}=j_{x}\displaystyle\frac{1+\varphi(x)}{\sqrt{(1+\varphi(x))^{2}-1-(a(x))^{2}}}, \;\;\;\varphi(0)=0, \;\;\varphi(1)=\varphi_{L},\\
\frac{d^{2}a}{dx^{2}}=j_{x}\displaystyle\frac{a(x)}{\sqrt{(1+\varphi(x))^{2}-1-(a(x))^{2}}},\;\;\;a(0)=0, \;\;\;a(1)=a_{L}.
\end{array} \label{eq:systemI}
\end{equation}

% 2
Where $j_x > 0,$  $\alpha \in [0,1], $  $\varphi$ is a potential of electric field and 
the potential of magnetic fields is $a.$

% 3

Let us now define the effective potential by 
$\theta(x) = (1+\varphi(x))^2 -1 - (a(x))^2.$

\section{Main Results}

% 4
We introduce the substitution $ u = 1 + \varphi, $ $a=v$ and reduce problem  to the following form:
\begin{equation} 
\begin{array}{l}
\label{31}
u^{\prime \prime} = j_x \frac{u}{\sqrt{u^2-1 - v^2}},  \, 
u(0) = 1, u(1) = \varphi_{L} + 1 = \alpha, u^{\prime}(0) = 0, \\
v^{\prime \prime} = j_x \frac{v}{\sqrt{u^2 - 1 - v^2}},
v(0)=0, v(1) = a_{L}, v^{\prime}(0) = \beta >0. 
\end{array}
\end{equation}

\begin{defi}
A function $(u,v) = (u(x),v(x))$ is a solution of the initial value problem (IVP) on $[0,\varepsilon)$ if:
\begin{itemize}
\item $u,v \in C^{1}[0,\varepsilon) \cap C^{2}(0,\varepsilon)$,
\item $\theta(u(x),v(x))>0$ for $x\in (0,\varepsilon)$,
\item $(u,v)$ satisfies the differential equations \eqref{31},
\item the initial conditions hold.
\end{itemize}
\end{defi}

\begin{prop}\label{prop:effectivePot}
Let $(u,v)$ be a solution of the IVP on $[0,\varepsilon)$ and define $\theta = u^{2} - 1 - v^{2}$. Then:
\begin{itemize}
\item $\theta(x) \in C^{1}[0,\varepsilon)\cap C^{2}(0,\varepsilon)$,
\item $\theta(0) = \theta'(0) = 0$ and $\theta(x)>0$ on $(0,\varepsilon)$,
\item $\theta$ satisfies the differential equation:
\begin{equation*}
\theta'' = j_{x}\left(6\sqrt{\theta} + \frac{2}{\sqrt{\theta}} - 4\gamma\right), \quad \gamma = -\frac{1}{2j_{x}}\beta^{2}.
\end{equation*}
\end{itemize}
\end{prop}

\begin{thm}
For every $\gamma\in \mathbb{R}$ there exists a unique solution on $[0,\infty)$ of the initial value problem:
\begin{equation}
D'' = j_{x}\left(6\sqrt{D} + \frac{2}{\sqrt{D}} - 4\gamma\right), \quad D(0) = D'(0) = 0. \label{DU}
\end{equation}
The solution has different properties depending on the value of $\gamma$:
\begin{description}
\item[Case 1:] $\gamma < 2$: $D\in C^{1}[0,\infty)\cap C^{2}(0,\infty)$, $D(\infty) = \infty$.
\item[Case 2:] $\gamma = 2$: $D \in C^{1}[0,\infty)\cap C^{2}(0,\infty)$, $D(\infty) = 1$.
\item[Case 3:] $\gamma > 2$: $D\in C^{1}[0,a]\cap C^{2}(0,a]$, can be extended periodically.
\end{description}
\end{thm}
Employing Banach's fixed-point theorem,  the existence of a unique solution to \eqref{DU} is proved.

\begin{prop}[Child-Langmuir Law]
Let $0< c\leq j_{x}\leq j_{x}^{max}$, $a=0$. Then equation 
\begin{equation*}
\varphi^{''}= j_{x}\frac{1+\varphi}{\sqrt{\varphi(2+\varphi)}}, \;\; \varphi(0)=0, \; \varphi(1) = \varphi_{L}
\end{equation*}
has a lower positive solution $u_{0}=\delta^{2}x^{4/3}$ if $4\delta^{3}\geq 9 j_{x}^{max}(1+\delta^{2}) / \sqrt{2+ \delta^{2}}$ and an upper positive solution $u^{0}= \alpha + \beta x$ with $\varphi_{L} \geq \delta^{2}$.
\end{prop}

\subsection*{Insulated case}

%Here we will define how we get the First order ODE from \rm (I)

For the study of the magnetically insulated diode (MID), $\theta < 0$ we want to analyze the conditions under which this phenomenon occurs.

Considering the system \eqref{eq:systemI} and the effective potential $\theta$ defined in Proposition \ref{prop:effectivePot}, we get the following first-order ODE after some simple algebraic and calculus-based manipulation:
\begin{equation*}
    (\theta')^2 = k \theta + 8 j_x \theta ^{\frac{3}{2}} + 8 j_x \theta ^{\frac{1}{2}} + 4\beta^2 \quad , \quad  \theta(0)=0,
\end{equation*}
where $k$ is an integration constant that we got during the order reduction process. As we want to know when the electrons start to deflect, we will use this expression to find where the $\theta$ values decrease. We will call this point $x_d$. Then

\begin{equation*}\label{eq:nonlinear}
    (\theta(x_d)')^2 = \theta_d' = k \theta_d + 8 j_x \theta_d ^{\frac{3}{2}} + 8 j_x \theta_d ^{\frac{1}{2}} + 4\beta^2 = 0,
\end{equation*}

\begin{prop}[$\theta$ and $u$ solutions relationship]\label{prop:rel}
    Let $\hat{k} = \frac{k}{8j_x}$ and $\hat{\beta} = \frac{\beta^2}{2j_x}$, and consider the equations:
    \begin{equation*}
        u^3 + \hat{k}u^2 + u + \hat{\beta}^2 = 0, \qquad \theta^{\frac{3}{2}} + \hat{k}\theta + \theta^{\frac{1}{2}} + \hat{\beta}^2 = 0.
    \end{equation*}
    Then, any solution $u$ of the first equation with $\operatorname{Re}(u) > 0$ induces a solution $\theta = u^2$ of the second equation.
\end{prop}

As solving the equation described in terms of $u$ induces the desired solutions for $\theta_d$, the efforts can be focused on solving the cubic expression, which is more tractable than the corresponding radical expression.

\begin{prop}\label{prop:u_sol} Let $\Delta_u = 18\hat{k}\hat{\beta} + \hat{k}^2 - 4 - 4 \hat{k}^3\hat{\beta} - 27\hat{\beta}^2$, then the cubic equation $u^3 + \hat{k}u^2 + u + \hat{\beta}^2 = 0$ has the following solutions in $\mathbb{C}$:

\begin{itemize}
    \item $\Delta_u < 0$:

    $$u_1 = -\frac{\hat{k}}{3} + \frac{\sqrt[3]{4}}{18} \left( \sqrt[3]{A_1 + A_2} + \sqrt[3]{A_1 - A_2}\right)$$
    $$u_2 = \left[-\frac{\hat{k}}{3} - \frac{\sqrt[3]{4}}{18} \left( \sqrt[3]{A_1 + A_2} +  \sqrt[3]{A_1 - A_2}\right)\right] + $$ $$+ \left[  \frac{\sqrt{3}\sqrt[3]{4}}{36} \left( \sqrt[3]{A_1 + A_2} - \sqrt[3]{A_1 - A_2}\right) \right]i$$
    $$u_3 = \left[-\frac{\hat{k}}{3} - \frac{\sqrt[3]{4}}{18} \left( \sqrt[3]{A_1 + A_2} + \sqrt[3]{A_1 - A_2}\right)\right] +$$ $$+ \left[ - \frac{\sqrt{3}\sqrt[3]{4}}{36} \left( \sqrt[3]{A_1 + A_2} - \sqrt[3]{A_1 - A_2}\right) \right]i$$

    where $A_1 = -54\hat{k}^3+243\hat{k}-729\hat{\beta}$ and $A_2 = \sqrt{A_1^2 + 2916(3-\hat{k}^2)^3}$ 

    \item $\Delta_u = 0$ , $\hat{k} = \pm \frac{\sqrt{3}}{9}$ , $\hat{\beta} = \pm \sqrt{3}$:

    $$u_1 = u_2 = u_3 = \mp \frac{\sqrt{3}}{3}$$

    \item $\Delta_u = 0$ , $\hat{k} \neq \pm \frac{\sqrt{3}}{9}$ , $\hat{\beta} \neq \pm \sqrt{3}$:

    $$u_1 = \frac{\hat{k}^3 - 4 \hat{k} + 9 \hat{\beta}}{3- \hat{k}^2}$$
    $$u_2 = u_3 = \frac{-\hat{k}+9\hat{\beta}}{ 2\hat{k}^2-6}$$

    \item $\Delta_u > 0$ , $\hat{k} \neq \pm \frac{\sqrt{3}}{9}$ , $\hat{\beta} \neq \pm \sqrt{3}$:

    $u_1 = A_3 \cos{\left( \frac{1}{3} \arccos(A_4) \right)} - \frac{\hat{k}}{3}, $
    $u_2 = A_3 \cos{\left( \frac{1}{3} \arccos(A_4) + \frac{2\pi}{3} \right)} - \frac{\hat{k}}{3},$
    $u_3 = A_3 \cos{\left( \frac{1}{3} \arccos(A_4) + \frac{4\pi}{3}\right)} - \frac{\hat{k}}{3},$
    where $A_3 = \frac{2}{3}\sqrt{\hat{k}^2-3}$ and $A_4 = \frac{4\hat{k}^3-9\hat{k}+27\hat{\beta}}{6-2\hat{k}^2}(\hat{k}^2-3)^{-1/2}$ 
    
\end{itemize}
\end{prop}

With the analytical description of the solutions of the cubic equation in Proposition \ref{prop:rel} and the method to induce the solutions for $\theta_d$,  the behavior of the solutions can be described by understanding and manipulating the parameters $\hat{k}$ and $\hat{\beta}$.

%As we define our first order ODE, we define the process to get the cubic equation and its roots

%We define the solutions in the space of theta, solutions written in terms of k hat and b hat

\section{Bifurcation diagrams}

In order to understand how the parameters $\hat{k}$ and $\hat{\beta}$ influence the existence or non-existence of the solutions required for the magnetic insulated diode, the different bifurcation diagrams over the solutions in the $\theta$ and $u$ space can be employed.

%Here goes how diagrams were built

\subsection*{Construction of Bifurcation Diagrams via Analytical Solutions}

The bifurcation diagrams presented are constructed by analytically solving the steady-state condition of the system, specifically the differential equation $\theta' = 0$. This approach focuses on identifying equilibrium points of the system as a function of a control parameter, such as magnetic field strength.

\subsubsection*{Methodology for Constructing the Bifurcation Diagram}

\begin{enumerate}
  \item {Analytical Determination of Equilibrium Points}:  
  The condition $\theta' = 0$ implies that the system is at equilibrium. By setting the derivative to zero, the resulting algebraic equation is solved analytically to find expressions for $\theta$ in terms of the control parameter. This process yields the equilibrium solutions of the system.

  \item {Parameter Variation and Solution Evaluation}:  
  The analytical expressions for $\theta$ are evaluated across a range of values for the control parameter. This step involves substituting different parameter values into the analytical solutions to compute the corresponding equilibrium points.

  \item {Plotting the Bifurcation Diagram}:  
  The computed equilibrium points are plotted against the varying control parameter using Python's \texttt{matplotlib} library. Each point on the diagram represents a steady-state solution for a specific parameter value. Connecting these points reveals the structure of the solution branches and highlights bifurcation points where qualitative changes in the system's behavior occur.
\end{enumerate}

This analytical approach to constructing bifurcation diagrams provides a clear and precise understanding of the system's steady-state behavior, facilitating the identification of critical parameter values that lead to qualitative changes in dynamics. The Python code used to generate the bifurcation diagrams is \href{https://github.com/omardtl24/Mag_Isol_Numerical}{available online}.

%End of the explanaition

\subsection*{2D diagrams for $u$ and $\theta$}

%Here we explore the bifurcation diagrams over the solutions in u and we explain the behaviour related to the cubic equation

In this first set of diagrams, one parameter was fixed at a specific value, while the other was varied over the interval $[-5,5]$. The fixed values were chosen to allow the existence of a real solution to the cubic equation with multiplicity three.

\begin{figure}[h]
    \centering
    \includegraphics[width=1.0\textwidth]{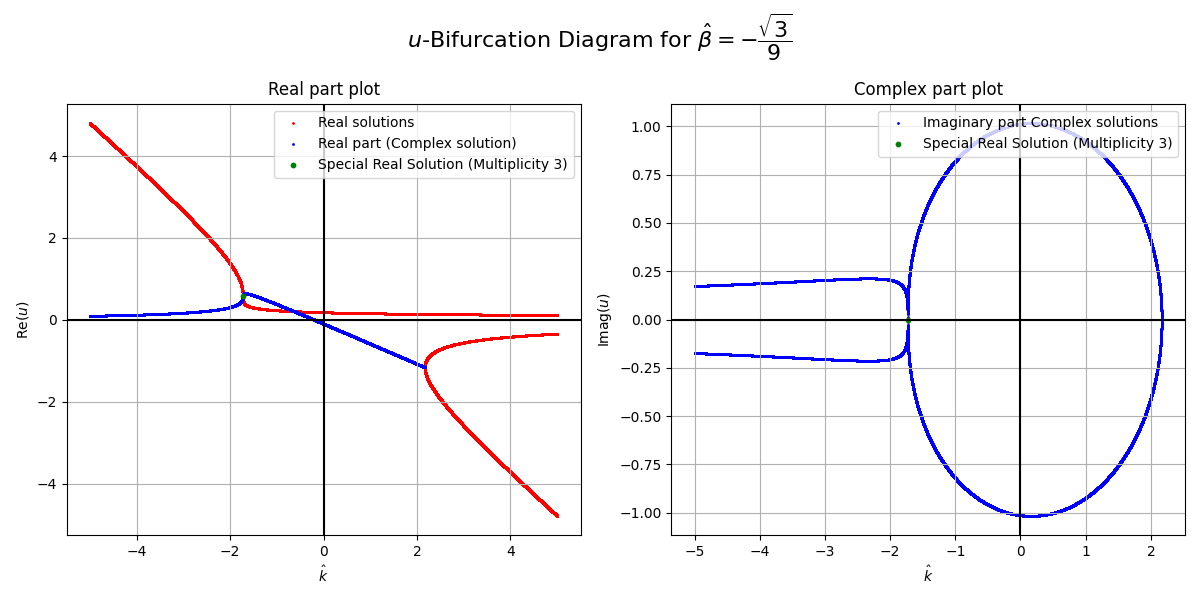}
    \caption{Bifurcation diagram of the solutions of the cubic equation related to $u$ centered on the 3-multiplicity real solution}
    \label{fig:u_2d_plot}
\end{figure}

In Figure \ref{fig:u_2d_plot}, we observe the behavior of the solutions in the $u$-space. As predicted by theory, the number of solutions in the complex plane ranges from one to three. However, since complex solutions occur in conjugate pairs, the plot reveals loop-like structures, which are of particular interest for further analysis.

%Here we explain the bifurcation diagrams over the solutions in theta, we show the interesting value 1/3 as it shows meaningful insights over the bifurcation process

\begin{figure}[h]
    \centering
    \includegraphics[width=1.0\textwidth]{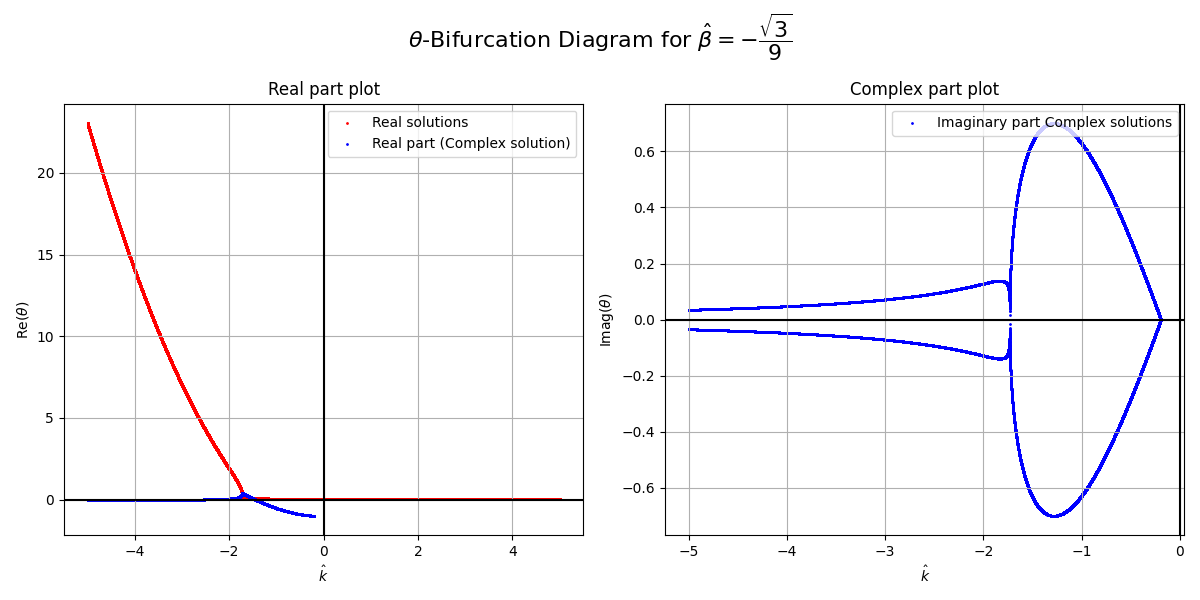}
    \caption{Bifurcation diagram of the solutions of the cubic equation related to $\theta$ centered on the 3-multiplicity real solution}
    \label{fig:theta_2d_plot}
\end{figure}

In Figure \ref{fig:theta_2d_plot}, we observe the behavior of the solutions in the $\theta$-space. In this new setting, some solutions are lost compared to the $u$-space due to the real part constraint defined in 
Proposition \ref{prop:rel}. Nevertheless, the loop-like behavior persists, which gives hope for its potential application to the magnetically insulated scenario.

\subsection*{3D surface diagrams for $u$ and $\theta$}

%Present the final surface regions gotten in 3D and explain the loop behaviour found on them.

The bifurcation diagrams in the 2D space revealed several interesting features. However, to obtain a more comprehensive understanding of the phenomenon, we now examine the behavior of the solutions in a 3D space. In this setting, both parameters are allowed to vary freely within the interval $[-5, 5]$, and the $z$-axis represents either the real or imaginary part of the solutions.

\begin{figure}[H]
    \centering
    \includegraphics[width=1.0\textwidth]{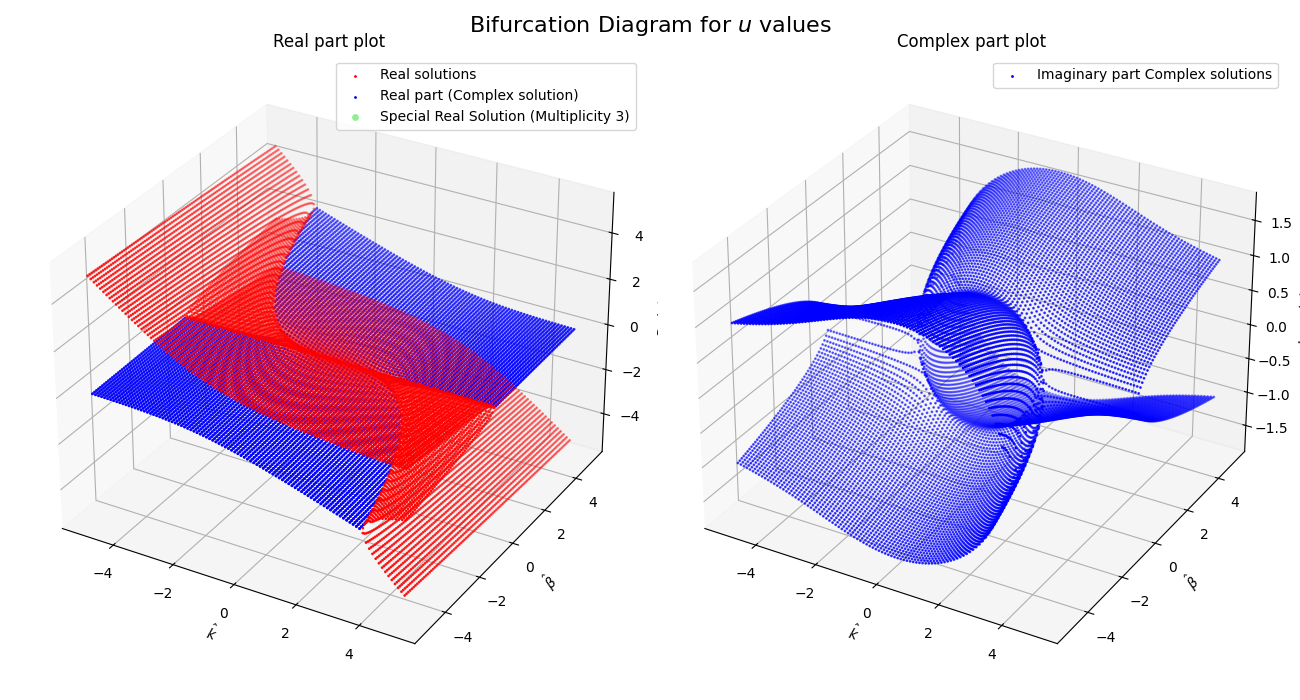}
    \caption{Surface gotten from the parameters $\hat{k}$ and $\hat{\beta}$ for the solutions of the cubic equation}
    \label{fig:u_3d_plot}
\end{figure}

In Figure \ref{fig:u_3d_plot}, we analyze the behavior of the solutions in the $u$-space. As a first observation, we note that for every pair of values  $(\hat{k}, \hat{\beta})$, there is at least one solution represented on the surface. Moreover, the loop-like behavior previously observed in the 2D case begins to manifest in a more structured and continuous form in the 3D surface, suggesting an underlying pattern that warrants further investigation.

\begin{figure}[H]
    \centering
    \includegraphics[width=1.0\textwidth]{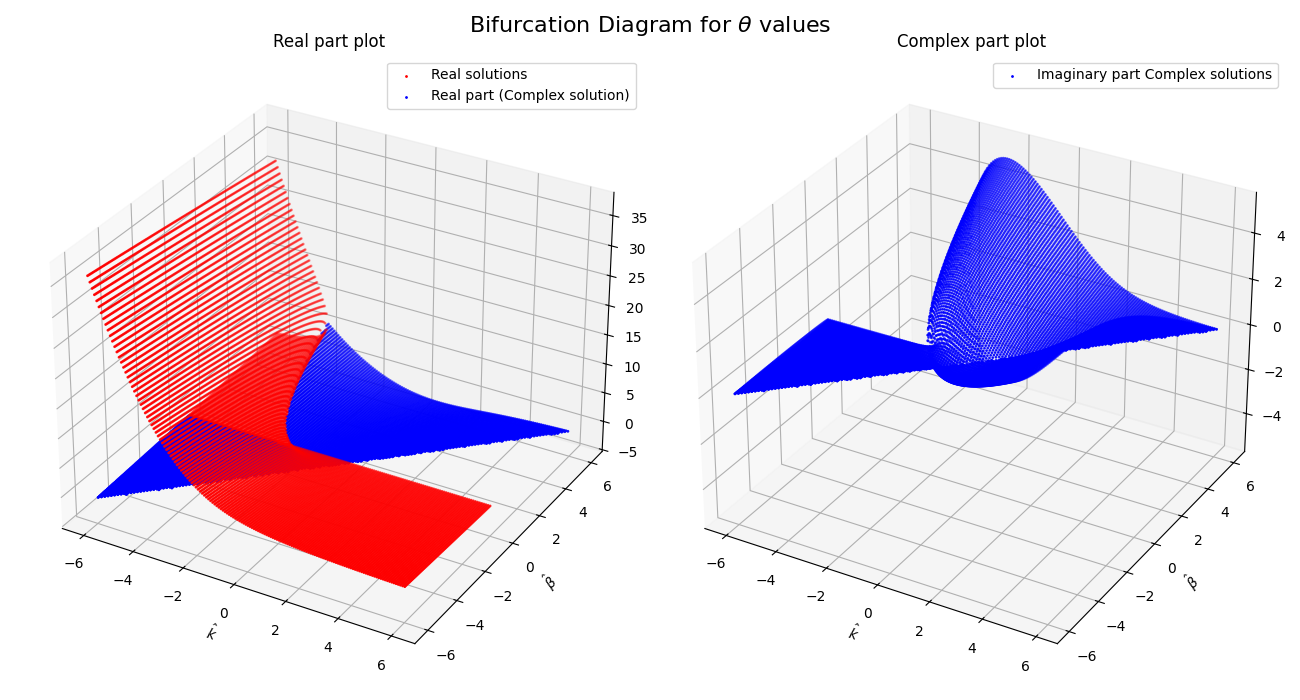}
    \caption{Surface gotten from the parameters $\hat{k}$ and $\hat{\beta}$ for the solutions of the stability of the $\theta$-based ODE}
    \label{fig:theta_3d_plot}
\end{figure}

In Figure \ref{fig:theta_3d_plot}, we analyze the behavior of the solutions in the $\theta$-space. Notably, we observe that the solutions with a negative real component are trimmed due to the constraint defined in Proposition \ref{prop:rel}. Additionally, there are empty regions that exhibit a likely linear pattern, which could be an interesting subject for further investigation.

\section{Conclusion}

This study yielded critical insights into the boundary problem governing MID through a rigorous analysis of the singularly perturbed Vlasov-Maxwell system. By reducing the system to a nonlinear ordinary differential equation framework,  the existence of physically admissible solutions under well-defined parametric constraints is proved. The bifurcation analysis underscored the important role of the free boundary point in shaping the solution topology, revealing distinct regimes dependent on system parameters. Numerical simulations further elucidated the mechanistic interplay between magnetic insulation and electron dynamics, particularly the redirection of electron trajectories toward the cathode under varying insulation thresholds. These findings emphasize the need for continued exploration of parameter-driven bifurcation phenomena, which may advance both theoretical and applied understanding of MID in high-energy physics and related fields. The bifircation analysis enables robust design strategies to avoid unstable regimes that could 
compromise device longevity and efficiency.

This study offers significant contributions to the design of high-performance power converters through the systematic optimization of MID operational parameters. Specifically, the analysis provides three key engineering advantages:  

\begin{enumerate}
\item Enhanced Operational Stability: By identifying the optimal working point  of the magnetically insulated diode, the proposed methodology ensures robust converter stability under dynamic load conditions. This equilibrium point minimizes oscillatory behavior and mitigates instability risks inherent in high-power systems.  

\item Improved Energy Density: The approach enables the accurate calibration of the diodes rated power capacity, directly enhancing energy density. This optimization aligns the converters output with application-specific requirements, thereby maximizing efficiency while avoiding overdesign.  

\item Cost-Effective Hardware Implementation: Traditional designs often rely on conservative maximum capacity estimates, leading to oversized components and inflated costs. In contrast, the derived optimal power threshold  allows for tailored diode specifications, reducing material expenditures without compromising performance or reliability.  

\end{enumerate}
Collectively, these insights advance the practical deployment of magnetically insulated diodes in next-generation converters.

\section*{Acknowledgments}
The study was done with a support of the state assignment No~FZZS-2024-0003 of Ministry of Science and Higher Education of Russian Federation.
The authors are grateful to Jurgen Batt, Pierre Degond and Walter Strauss for stimulating discussions.

%\section{Reference List} We recommend using the following samples for references. The list of references should be in alphabetic order. If available, please always include DOIs as full DOI links in your reference list (e.g. �https://doi.org/abc�).


\begin{thebibliography}{999}



\bibitem{2}
Ben Abdallah N.,  Degond P.,   Mehats F. { Mathematical Models of Magnetic Insulation}. \textit{Rapport interne N 97.20}. 1997, MIP, Universite Paul Sabatier, Toulouse, France.

\bibitem{3}
Degond P.,  Raviart P.-A., { An asymptotic analysis of the one-dimensional Vlasov-Poisson system: the Child-Langmuir law},  \textit{Asymptotic Anal.}, 1991, vol. 4, no. 3, pp. 187--214.

\bibitem{1}
Langmuir I., Compton K.T., { Electrical discharges in gases. Part II. Fundamental phenomena in electrical discharges}, \textit{ Rev. Mod. Phys.}, 1931,  vol. 3, no. 2,  pp. 191--257.

\bibitem{sym20}
Rojas E.M., Sidorov N.A., Sinitsyn A.V.  A boundary value problem for noninsulated magnetic regime in a vacuum diode, \textit{Symmetry},  2020, vol. 12, no. 4,  617.

\bibitem{sidsid}
Sidorov N.A., Sidorov D.N. Small solutions of nonlinear differential equations near branching points. 
\textit{Russ Math.}, 2011, vol. 55, pp. 43--50. "https://doi.org/10.3103/S1066369X11050070"

\bibitem{sidsidsin}
Sidorov N.,  Sidorov D., Sinitsyn A. Toward General Theory of Differential-Operator and Kinetic Models,  Book Series: \textit{World Scientific Series on Nonlinear Science Series A}, vol. 97, eds. Prof. L. Chua, World Scientific, Singapore, London, 2020 , 496 pp.

\bibitem{7}
Sinitsyn A.V. {Positive solutions of a nonlinear singular boundary value problem of magnetic insulation}, \textit{Mathematical modelling}, 2001, vol. 13, no. 5,  pp. 37--52.

\bibitem{8} 
Sidorov D., Rojas E., Sinitsyn A., Sidorov D.  Approximation and Regularisation Methods for Operator-Functional Equations, Book Series: \textit{Advances in Mathematics for Applied Sciences}, vol. 95, World Scientific, Singapore, London, 2025 , 248 pp.


%\bibitem{Kr1965} Krni\'c L. Types of Bases in the Algebra of Logic.  \textit{Glasnik Matematicko-Fizicki i Astronomski}, ser 2, 1965, vol. 20, pp. 23-32.
%
%
%\bibitem{L2008}Lau  D.,  Miyakawa M. Classification and enumerations of bases in $P_k(2)$. \textit{Asian-European Journal of Mathematics}, 2008, vol. 1, no. 2, pp. 255-282.
%
%
%
%\bibitem{M1990}  Miyakawa M., Rosenberg I., Stojmenovi\'c I. Classification of Three-valued logical functions preserving 0. \textit{Discrete Applied Mathematics}, 1990,  vol. 28, pp. 231-249. https://doi.org/10.1016/0166-218X(90)90005-W




\end{thebibliography}
\end{document}